\documentclass[12pt]{article}
\usepackage{latexsym}
\usepackage{amssymb,amsmath,amsthm}
\usepackage{addfont} 
\addfont{OT1}{rsfs10}{\rsfs}
\usepackage{graphicx}
\RequirePackage{tikz}
\usepackage{epsfig}
\usepackage{multicol}
\usepackage{pst-grad} 
\usepackage{pst-plot} 
\usepackage{pst-node} 
\usepackage{pst-tree}
\usetikzlibrary{shapes, arrows}
\usepackage[all]{xy}
\graphicspath{{fig/}}
\usepackage{subfig}
\usepackage{amsfonts}

\theoremstyle{definition}
\newtheorem{theorem}{Theorem}[section]

\newtheorem{corollary}[theorem]{Corollary}

\newcommand{\T}{\mathrm}

\begin{document}

\title{On the comparison of the distinguishing coloring and the locating coloring of graphs}

\author{M. Korivand $^{1}$, A. Erfanian  $^{1}$, and  Edy T.  Baskoro $^{2}$  }

\date{}

\maketitle

\begin{center}

$^{1}$ Department of Pure Mathematics, Ferdowsi University of Mashhad, \\
P.O.Box 1159-91775, Mashhad, Iran. \\
e-mail: {\tt korivand@mail.um.ac.ir }, {\tt erfanian@math.um.ac.ir}

$^{2}$ Combinatorial Mathematics Research Group, ITB, Jl. Ganesha 10, Bandung. \\
e-mail: {\tt ebaskoro@math.itb.ac.id}

\end{center}

\begin{abstract}
Let $G$ be a simple connected graph. Then $\chi_L(G)$ and $\chi_{D}(G)$ will denote the locating chromatic number and the distinguishing chromatic number  of $G$, respectively. In this paper, we investigate a comparison between  $\chi_L(G)$ and $\chi_{D}(G)$. In fact, we prove that  $\chi_{D}(G)\leq \chi_L(G)$. Moreover, we
determine some types of graphs whose locating and distinguishing chromatic
numbers are equal. Specially, we characteristic all graph 
$G$ 
with property that 
$ \chi_{D}(G) = \chi_{L}(G) = 3 $.
\end{abstract}


\medskip\noindent
{\bf AMS Subj.\ Class:} 05C15.

\section{Introduction} 
One of the important and applicable concepts in graph theory is graph coloring. The subject graph coloring is one of the best known, popular, and extensively researched subjects in the field of graph theory, having many conjectures, which are still open and studied by various mathematicians and computer scientists along with the world. In this article, we will discuss two types of graph coloring that we remind them as following.

Let $G$ be a graph and $c$ be a proper $k$-coloring of graph $G$ and $\pi=(V_1, V_2, \cdots, V_k)$ where $V_i$ is the set of all vertices colored by $i$ for $1\leqslant i\leqslant k$. The {\it color code of vertex $v$} with respect to $\pi$, denoted by $c_{\pi}(v)$, is defined as the ordered $k$-tuple $(\T{d}(v, V_1), \T{d}(v, V_2), \cdots, \T{d}(v, V_k))$ such that $\T{d}(v, V_i)$ is the minimum distance from $v$ to each  vertex in $V_i$ for $1\leqslant i\leqslant k$. If all  vertices of $G$ have distinct color codes, then $c$ is called  a {\it locating coloring of $G$}. The locating chromatic number, $\chi_L(G)$, is the minimum number of colors needed in a locating coloring of $G$. 

The concept of locating coloring was first introduced by Chartrand et. al.  \cite{chart} in 2002. Since then, the locating coloring number has been the subject of many researchers, for more details see \cite{tri, erwin, beh}. 

A $k$-coloring of a graph $G$ is said  to be a {\it distinguishing $k$-coloring of $G$} if it is the proper $k$-coloring of $G$ and the identity automorphism is the only color-preserving automorphism of $G$. A distinguishing chromatic number $\chi_{D}(G)$ of $G$ is the least $k$ such that $G$ has a distinguishing $k$-coloring.

In 2006, Collins and Trenk \cite{trenk} introduced  the distinguishing chromatic number of a graph. Many authors  obtained more results on the distinguishing chromatic number and related subjects, see \cite{choi, fisher}.

The locating coloring of a graph is dealing with the distance between vertices of a graph. The distinguishing coloring is discussing about vertices and automorphisms which are distance-preserving. So, it would be interested to give a relation between the above two kinds of coloring. In this article, we will compare the locating and distinguishing colorings. In fact, we prove that any locating coloring of a connected graph is a distinguishing coloring. Note that in some cases, for example complete multipartite graphs, the above two colorings are the same, but there are many cases of a graph  $G$ that $\chi_{D}(G) < \chi_L(G)$. It is the most of our interest to see that when $\chi_{D}(G)=\chi_L(G)$ or $\chi_{D}(G)<\chi_L(G)$. In this paper, we state some results for sufficient conditions that a distinguishing coloring is a locating coloring. Also, 
we characterize all graph 
$ G $ 
with order 
$ n $
such that 
$ \chi_{D}(G) = \chi_{L}(G) = 3 $.

All graphs mentioned in this paper are assumed to be simple, undirected, and connected. So, we do not need to state these conditions in our results. Moreover,  we write $N(v)$ for the set of neighbours of a vertex $v$ in a graph and other notations and terminologies are standard and one can refer to \cite{bondy}.

\label{sec:intro}

\section{Main Results} 

In this section, we are going to compare the locating chromatic number and the distinguishing chromatic number of  a graph. First of all,  let us state the following two theorems from \cite{chart, trenk} which states necessary and sufficient conditions for a graph $G$ such that $\chi_L(G)=|V(G)|$ and similarly whenever $\chi_{D}(G)=|V(G)|$. 
\begin{theorem}[\cite{chart}]\label{th2.1}
Let $G$ be a connected graph of order $n\geq 3$. Then $\chi_L(G)=n$ if and only if $G$ is a complete multipartite graph. 
\end{theorem}
\begin{theorem}[\cite{trenk}]\label{th2.2}
Let $G$ be a graph. Then $\chi_{D}(G)=|V(G)|$ if and only of $G$ is a complete multipartite graph.
\end{theorem}
Now,  we compare the locating and distinguishing coloring of a graph in the following theorem.
\begin{theorem}\label{th2.3}
Let $G$ be a graph. Then any locating coloring of $G$ is a distinguishing coloring of $G$.
\end{theorem}
\begin{proof}
Let $[\chi_L(G)]=\{V_1, V_2, \cdots, V_k\}$ be a color classes of a locating coloring of $G$. If  $|V_i|=1$ for all $i$, $1\leq i\leq k$, then  $\chi_L(G)=|V(G)|$ by Theorem \ref{th2.1}. Hence, $G$ is a complete multipartite graph. So, by Theorem \ref{th2.2} $\chi_{D}(G)=|V(G)|$  and $[\chi_L(G)]$ is a distinguishing coloring of $G$.

Thus, let $|V_{j}|\geq 2$ for some $j$, $1\leq j\leq k$. We claim that $[\chi_L(G)]$ is a distinguishing coloring of $G$. Suppose, to the contrary, that $[\chi_L(G)]$ is not a distinguishing coloring of $G$.  Hence, there exists a non-identity automorphism $f$  that preserves the color classes. Without loss of generality, we may assume that $V_j=\{a_1, a_2, \cdots, a_m\}$ with $m\geq 2$ and $f(a_1)=a_2$. Let for each $i$, $1\leq i\leq k$,  $\T{d}(a_1, V_i)=\T{d}(a_1, b_i)$ for some $b_i\in V_i$. Then
\begin{equation}\label{eq1}
\T{d}(a_1, V_i)=\T{d}(a_1, b_i)= \T{d}(f(a_1), f(b_i))=\T{d}(a_2, f(b_i)).
\end{equation}
Let $\T{d}(a_2, V_i)=\T{d}(a_2, c_i)$, where $ c_i\in V_i$. 

First, assume that $ f(b_i) \neq c_i $, we have 
\begin{equation}\label{eq2}
\T{d}(a_2, V_i)=\T{d}(a_2, c_i)=\T{d}(f^{-1}(a_2), f^{-1}(c_i))=\T{d}(a_1, f^{-1}(c_i)).
\end{equation}
Since $f$ preserves all color classes, it follows that $f(b_i)\in V_i$ and $f^{-1}(c_i)\in V_i$. Hence, $\T{d}(a_1, f^{-1}(c_i))\geq \T{d}(a_1, b_i)=\T{d}(a_1, V_i)$ and  $\T{d}(a_2, f(b_i))\geq \T{d}(a_2, c_i)=\T{d}(a_2, V_i)$. Then from \eqref{eq1} and \eqref{eq2} we have $\T{d}(a_1, V_i)=\T{d}(a_2, V_i)$ for all $1\leq i\leq k$. It means that the color codes of $a_1$ and $a_2$ are the same and it is a contradiction. 

Finally, suppose that   $ f(b_i) = c_i $,  we have $ \T{d}(a_1, f^{-1}(c_i)) = \T{d}(a_1, b_i) $, this implies that $\T{d}(a_1, V_i)=\T{d}(a_2, V_i)$ for all $1\leq i\leq k$, and similar to above we reach a contradiction.
 Hence, the proof is completed.
\end{proof}
By Theorem \ref{th2.3}, the following result is obtained, directly.
\begin{corollary}\label{cor2.5}
Let $G$ be a graph. Then $\chi_{D}(G)\leq \chi_L(G)$.
\end{corollary}
We remind the metric dimension of a graph $G$. For an ordered set $W=\{ w_1, w_2, \cdots, w_k\}$ of vertices in a connected graph $G$ and a vertex $v$ of $G$, the ordered $k$-tuple $r_W(v)$ of $v$ with respect to $W$ is defined by  
\[r_W(v)=(\T{d}(v, w_1), \T{d}(v, w_2), \cdots, \T{d}(v, w_k)),\]
where $\T{d}(v, w_i)$ is the distance between $v$ and $w_i$, $1\leq i\leq k$. The set $W$ is a resolving set for $G$ if the $k$-tuples $r_W(v)$, $v\in V(G)$, are distinct. The {\it metric dimension of $G$} is the minimum cardinality of a resolving set for $G$ and is denoted by $\text{dim}(G)$.  These concepts were introduced independently in \cite{harary, slater}.  

 In \cite{chart}, Chartrand et. al. gave an upper bound for $\chi_L(G)$ in terms of $\chi(G)$ and $\text{dim}(G)$. They prove that $\chi_L(G)\leq \chi(G)+\text{dim}(G)$. It is obvious that the above upper bound is sharp. For instance, if $P_{2k+1}$ is a path of length $2k$, then $\chi_L(P_{2k+1})=3$, $\chi(P_{2k+1})=2$ and $\text{dim}(P_{2k+1})=1$. Now, by Theorem \ref{th2.3}, we can easily get the above upper bound for $\chi_{D}(G)$.
\begin{corollary}
For any connected graph $G$, $\chi_{D}(G) \leq \chi(G)+\text{dim}(G)$.
\end{corollary}
One can see that any upper bound for $\chi_L(G)$ will be an upper bound for $\chi_{D}(G)$ as well, similarly, any lower bound for $\chi_{D}(G)$ will be a lower bound for $\chi_L(G)$. The following  corollary is coming from this point of view. We can refer to \cite{chart}.
\begin{corollary}
If $G$ is a graph of order $n\geq 3$ with diameter  $\text{diam}(G)\geq 2$, then $\chi_{D}(G)\leq n-\text{diam}(G)+2$.
\end{corollary} 
Notice that any distinguishing coloring of $G$ is not necessarily to be a locating coloring of $G$. For example, Let $[\chi_{D}(P_7)]=\{ \{a_1, a_5\}, \{a_2, a_4, a_6\}, \{a_3, a_7\}\}$ be a distinguishing  coloring classes of $P_7$.  This coloring of $P_7$ is a distinguishing but is not locating coloring, because  $c_{\pi}(a_2)=c_{\pi}(a_4)=(1, 0, 1)$. Also, it can be seen that $\chi_{D}(P_7)=\chi_L(P_7)=3$.
\begin{figure}[h!]
\centering
\begin{tikzpicture}[scale=.8, thick, vertex/.style={scale=.6, ball color=black, circle, text=white, minimum size=.2mm},
Ledge/.style={to path={
.. controls +(45:2) and +(135:2) .. (\tikztotarget) \tikztonodes}}
]
\node [ label={[label distance=.1cm]-90:$a_1$}] [ label={[label distance=.1cm]90:$1$}]  (1) at (-6,0) [vertex] {};
\node [ label={[label distance=.1cm]-90:$a_2$}] [ label={[label distance=.1cm]90:$2$}] (2) at (-4,0) [vertex] {};
\node [ label={[label distance=.1cm]-90:$a_3$}] [ label={[label distance=.1cm]90:$3$}] (3) at (-2,0) [vertex] {};
\node [ label={[label distance=.1cm]-90:$a_4$}] [ label={[label distance=.1cm]90:$2$}] (4) at (0,0) [vertex] {};
\node [ label={[label distance=.1cm]-90:$a_5$}] [ label={[label distance=.1cm]90:$1$}]  (5) at (2,0) [vertex] {};
\node [ label={[label distance=.1cm]-90:$a_6$}] [ label={[label distance=.1cm]90:$2$}] (6) at (4,0) [vertex] {};
\node [ label={[label distance=.1cm]-90:$a_7$}] [ label={[label distance=.1cm]90:$3$}] (7) at (6,0) [vertex] {};
\draw (1) -- (2);
\draw (2) -- (3);
\draw (3) -- (4);
\draw (4) -- (5);
\draw (5) -- (6);
\draw (6) -- (7);
\end{tikzpicture}
\caption{A distinguishing coloring in $P_7$.}
\end{figure}
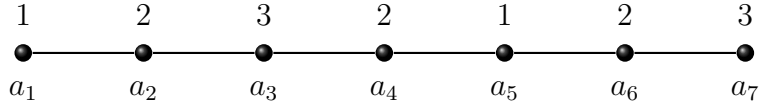

It is interesting to see that the difference between $\chi_{D}(G)$ and $\chi_L(G)$. It seems that the difference can be small or large and each one independent to another one. In other words, we can always find a graph $G$ with $\chi_{D}(G)=n$ but $\chi_L(G)=m$, where $n\leq m\leq 2n-1$. The following theorem states this fact. 
\begin{theorem}\label{f1}
There exists a graph $G$ having locating chromatic number $m$ and distinguishing chromatic number $n$, for all $2\leq n\leq m\leq 2n-1$.
\end{theorem}
\begin{proof}
If $ n = 2 $, then the result is travail. So, assume that $n>2$ is fixed. We consider the graph $G$ in Fig. 2 such that the degree of vertex 
$a$ is $(n-1)(n-1)$. It is easy to see that $\chi_L(G)=n$ and $\chi_D(G)=n$. Consider a locating $n$-coloring for $ G$. Without loss generality, we can give color 1 to vertex $a$. 
Hence, for every locating color class $V_i$, with color $i, 2 \leq i \leq n$, we have 
$|N(a) \cap V_i|=n-1$. So there exists exactly one vertex, calling $a_i$, in each $V_i$, $2 \leq i \leq n$, where it is adjacent to a pendant with color 1. Now, let $m$ be an integer such that $n \leq m \leq 2n-1$. 
Construct a new graph $G^*$ from $G$ by adding $m-n$ pendants to vertex $a$. Note that if $m=n$ then $G^*=G$. By coloring 
$m-n \; (>=1)$ new vertices with $2, 3, \cdots, m-n+1$ and preserving the previous coloring for other vertices, we obtain a minimum distinguishing coloring of $G^*$, and so $\chi_D(G^*)=n$. 
However, to get the right value of $\chi_L(G^*)$ consider any pendant vertex $b$ in $N(a)$. If it is colored by $i$, $1\leq i \leq n$, then $c_\pi(b)=c_\pi(a_i)$. Hence, we need $m-n$ new colors for coloring all pendant vertices in $N(a)$, and this implies that $\chi_L(G^*)=m$.

\begin{figure}[h!]
\centering
\begin{tikzpicture}[scale=.6 ,thick, vertex/.style={scale=.6, ball color=black, circle, text=white, minimum size=.2mm},
Ledge/.style={to path={
.. controls +(45:2) and +(135:2) .. (\tikztotarget) \tikztonodes}}
]
\node [ label={[label distance=.1cm]-90:$a$}]  (1) at (0,0) [vertex] {};
\node   (2) at (-1.6,1) [vertex] {};
\node   (3) at (-3.2, 2) [vertex] {};
\node   (4) at (-2,3.4) [vertex] {};
\node   (5) at (-1,1.7) [vertex] {};
\node   (6) at (1.8,.7) [vertex] {};
\node   (7) at (3.7,1.4) [vertex] {};
\node   (8) at (0.1,1.8) [vertex] {};
\node   (9) at(0.2 ,3.6)[vertex] {};
\node   [ label={[label distance=.2cm]-80:$.$}](10) at(1.2 ,3.8)[] {};
\node   [ label={[label distance=.2cm]-80:$.$}](11) at(1.6 ,3.5)[] {};
\node   [ label={[label distance=.2cm]-80:$.$}](12) at(2 ,3.2)[] {};
\draw (1) -- (2);
\draw (2) -- (3);
\draw (1) -- (5);
\draw (4) -- (5);
\draw (1) -- (6);
\draw (6) -- (7);
\draw (1) --(8);
\draw (8)--(9);
\end{tikzpicture}
\caption{}\label{fig2}
\end{figure}
\end{proof}

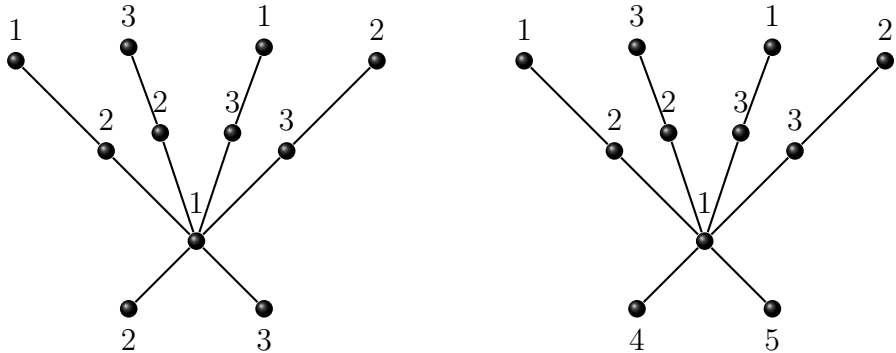
\begin{figure}[h!] \hspace{12mm}
\subfloat{\begin{tikzpicture}[scale=.6, thick, vertex/.style={scale=.6, ball color=black, circle, text=white, minimum size=.2mm},
Ledge/.style={to path={
.. controls +(45:2) and +(135:2) .. (\tikztotarget) \tikztonodes}}
]
 \node[label={[label distance=.1cm]90:$1$‌}] (a1) at (0, 0) [vertex] {};‌
\node[label={[label distance=.001cm]90:$2$‌}] (a2) at (-2, 2) [vertex] {};‌
\node[label={[label distance=.001cm]90:$1$}] (a3) at (-4, 4) [vertex] {};‌
\node[label={[label distance=.001cm]90:$2$}] (a4) at (-.80, 2.4) [vertex] {};‌
\node[label={[label distance=.001cm]90:$3$}] (a5) at (-1.5, 4.3) [vertex] {};‌
\node[label={[label distance=.001cm]90:$3$}] (a6) at (.80, 2.4) [vertex] {};‌
\node[label={[label distance=.001cm]90:$1$}] (a7) at (1.5, 4.3) [vertex] {}; ‌
\node[label={[label distance=.001cm]90:$3$}] (a8) at (2, 2) [vertex] {};‌ 
\node[label={[label distance=.001cm]90:$2$}] (a9) at (4, 4) [vertex] {};‌ 
\node[label={[label distance=.001cm]-90:$2$‌}] (b1) at (-1.5, -1.5) [vertex] {};‌
\node[label={[label distance=.001cm]-90:$3$}] (b2) at (1.5, -1.5) [vertex] {};‌
\draw(a1) to (a2) to (a3) ;
\draw(a1) to (a4) to (a5) ;
\draw(a1) to (a6) to (a7) ;
\draw(a1) to (a8) to (a9) ; 
\draw(b1) to (a1) to (b2) ;
\end{tikzpicture}}
\subfloat{
\hspace{12mm}
\begin{tikzpicture}[scale=.6 ,thick, vertex/.style={scale=.6, ball color=black, circle, text=white, minimum size=.2mm},
Ledge/.style={to path={
.. controls +(45:2) and +(135:2) .. (\tikztotarget) \tikztonodes}}
]
 \node[label={[label distance=.1cm]90:$1$‌}] (a1) at (0, 0) [vertex] {};‌
\node[label={[label distance=.001cm]90:$2$‌}] (a2) at (-2, 2) [vertex] {};‌
\node[label={[label distance=.001cm]90:$1$}] (a3) at (-4, 4) [vertex] {};‌
\node[label={[label distance=.001cm]90:$2$}] (a4) at (-.80, 2.4) [vertex] {};‌
\node[label={[label distance=.001cm]90:$3$}] (a5) at (-1.5, 4.3) [vertex] {};‌
\node[label={[label distance=.001cm]90:$3$}] (a6) at (.80, 2.4) [vertex] {};‌
\node[label={[label distance=.001cm]90:$1$}] (a7) at (1.5, 4.3) [vertex] {}; ‌
\node[label={[label distance=.001cm]90:$3$}] (a8) at (2, 2) [vertex] {};‌ 
\node[label={[label distance=.001cm]90:$2$}] (a9) at (4, 4) [vertex] {};‌ 
\node[label={[label distance=.001cm]-90:$4$‌}] (b1) at (-1.5, -1.5) [vertex] {};‌
\node[label={[label distance=.001cm]-90:$5$}] (b2) at (1.5, -1.5) [vertex] {};‌
\draw(a1) to (a2) to (a3) ;
\draw(a1) to (a4) to (a5) ;
\draw(a1) to (a6) to (a7) ;
\draw(a1) to (a8) to (a9) ; 
\draw(b1) to (a1) to (b2) ;

\end{tikzpicture}
}
\begin{center}
\caption{\ A distinguishing $ 3$-coloring (left) and locating $ 5$-coloring (right) of graph $ G$.}\label{fig3}
\end{center}
\end{figure}

In Figure \ref{fig3}, we give an illustration of Theorem \ref{f1}, when 
$ m = 5 $ 
and
$ n = 3$, 
for a graph 
$ G $ with $ \chi_{D}(G) = 3 $ and $ \chi_{L}(G) = 5$. 

In the following two theorems, we give some conditions for a graph $G$ such that the distinguishing chromatic number and locating chromatic number are equal.

\begin{theorem}
Let $G$ be a graph with $\chi_{D}(G)=|V(G)|-1$. Then any distinguishing coloring of $G$ is a locating coloring of $G$. 
\end{theorem}

\begin{proof}
Let $c$ be a distinguishing $\chi_D(G)$-coloring which induces the partition $\pi=(\{v_1\}, \{v_2\}, \cdots, \{v_{n-2}\}, \{v_{n-1},v_n\})$. The only vertices in $G$ with the same color are $v_{n-1}$ and $v_n$.
Therefore, it suffices to show that the color codes of these two vertices are distinct. Of course, $v_{n-1}$ and $v_n$ are not adjacent. Let $A=N(v_{n-1}) - N(v_{n})$ and $B=N(v_{n}) - N(v_{n-1})$.
Since under coloring $c$ there is no non-identity automorphism preserving the above color class, then
at least one of $A$ or $B$ is not empty, say $|A| \geq1$ and $v_i \in A$ for some $i$. 
Therefore, the $i^{th}$ component of $c_\pi(v_{n-1})$ is 1, but the one of $c_\pi(v_{n})$ is $\geq 2$. This means that the color codes of these two vertices are different.  Thus, $c$ is a locating coloring of $G$. 
\end{proof} 
\begin{corollary}
Every graph $G$ with $\chi_{D}(G) = |V (G)|-1$ has the locating chromatic number $|V (G)| - 1$. 
\end{corollary} 

\begin{theorem}
Let $G$ be a graph with $\chi_{D}(G) = |V (G)|-2$. Let $\pi$ be the partition induced by a distinguishing $\chi_{D}(G)$-coloring $c$ that satisfies one of the following conditions: 
\begin{itemize}
\item[(i)]
There is a color class of three vertices in $\pi$, or  
\item[(ii)] 
The two color $2$-classes $\{a, b\}$ and $\{c, d\}$ in $ \pi $ satistify that: 

$N_{S}(a) \neq N_{S}(b)$ and $N_{S}(c) \neq N_{S}(d)$ where $S$ is the set of all vertices in
singleton classes in $ \pi $. 
\end{itemize}
Then, $c$ is a locating coloring of $G$. 
\end{theorem}

\begin{proof}
Let $c$ be a distinguishing $\chi_D(G)$-coloring and $\pi$ be the induced partition by $c$.
Since  $\chi_{D}(G)=|V(G)|-2$, then there are two kinds of color classes namely, 
$\pi=(\{v_1\}, \{v_2\}, \cdots, \{v_{n-3}\}, \{v_{n-2}, v_{n-1}, v_n\})$ or
$\pi=(\{v_1\}, \{v_2\}, \cdots, \{v_{n-3}, v_{n-2}\}, \{v_{n-1}, v_n\})$.
In order to prove that $c$ is a locating coloring (if $c$ satisfies the required conditions), it suffices to show that the color codes of all vertices
in any non-singleton color class are distinct. 

Let $x$ and $y$ be two distinct vertices in a non-singleton color class $W$ in $\pi$. 
Let $A = N(x) - N(y)$ and $B = N(y) - N(x)$.
Since under coloring $c$ there is no non-identity automorphism preserving the above color class, then at least one of 
$A$ or $B$ is not empty, say $|A| \geq 1$. If $W = \{ v_{n-2}, v_{n-1}, v_n \} $ then there is $v_i$ in $A$ for some $ i$  
$ ( 1 \leq i \leq n - 3)$. 
Therefore, the $i^{th}$ component of $c_\pi(x)$ is 1, but the one of $c_\pi(y)$ is $\geq 2$. This means that the color codes of  
$ x $ 
and 
$ y $
are different. So, in this case, $c$ is also a locating coloring of $G$.

Now, let $W$ be either one of $\{v_{n-3}, v_{n-2} \}$ or $\{v_{n-1}, v_{n} \}$ and let $x, y \in W$. 
Let $S =\{v_1, v_2, \ldots, v_{n-4} \}$. 
Since $N_{S}(x) \neq N_{S}(y)$, then, there exists $v_i$ which is adjacent to exactly one of $\{x, y\}$, 
say $v_i \sim x$ but $v_i \nsim y$. Therefore, the $i^{th}$ component of $c_{\pi}(x)$ is $1$, but the one of $c_{\pi}(y)$ is $\geq 2$. This means that the color codes of $x$ and $y$ are distinct. So, in this case, $c$ is also a locating coloring of $G$.
\end{proof}

\section{Graphs with $ \chi_{D}(G)=\chi_{L}(G) = 3 $} 

In this section, we will further derive all graphs $G$ with $\chi_L(G)=\chi_D(G)$. In particular, we characterize all graphs $G$ with $\chi_L(G)=\chi_D(G)=3$. 
First, let $\mathcal{T}$ be a set of all trees $T$ with locating-chromatic number 3.
Baskoro and Asmiati (2013) characterized all trees on $n$ vertices ($n\geq 3$) with locating-chromatic number 3 as follows.\\

\begin{theorem}\cite{tri}
\label{Th-tree3}
A tree $T$ is in $\mathcal{T}$ if and only if $T$ is any subtree of one of the trees (A), (B) or (C) in Figure 4 containing vertices $X$, $Y$ and $Z$,
with $a \geq 0,b \geq 0,c \geq0, d\geq 0,e\geq 0,f\geq 0,i \geq 0,j \geq0,k\geq0,p\geq0$;  $g\geq 1, h \geq 1$; $i=j$ and $k=p$. \\
\end{theorem}

\begin{figure}[!h]
\begin{center}
\label{FigE1}
  \includegraphics[width=10cm]{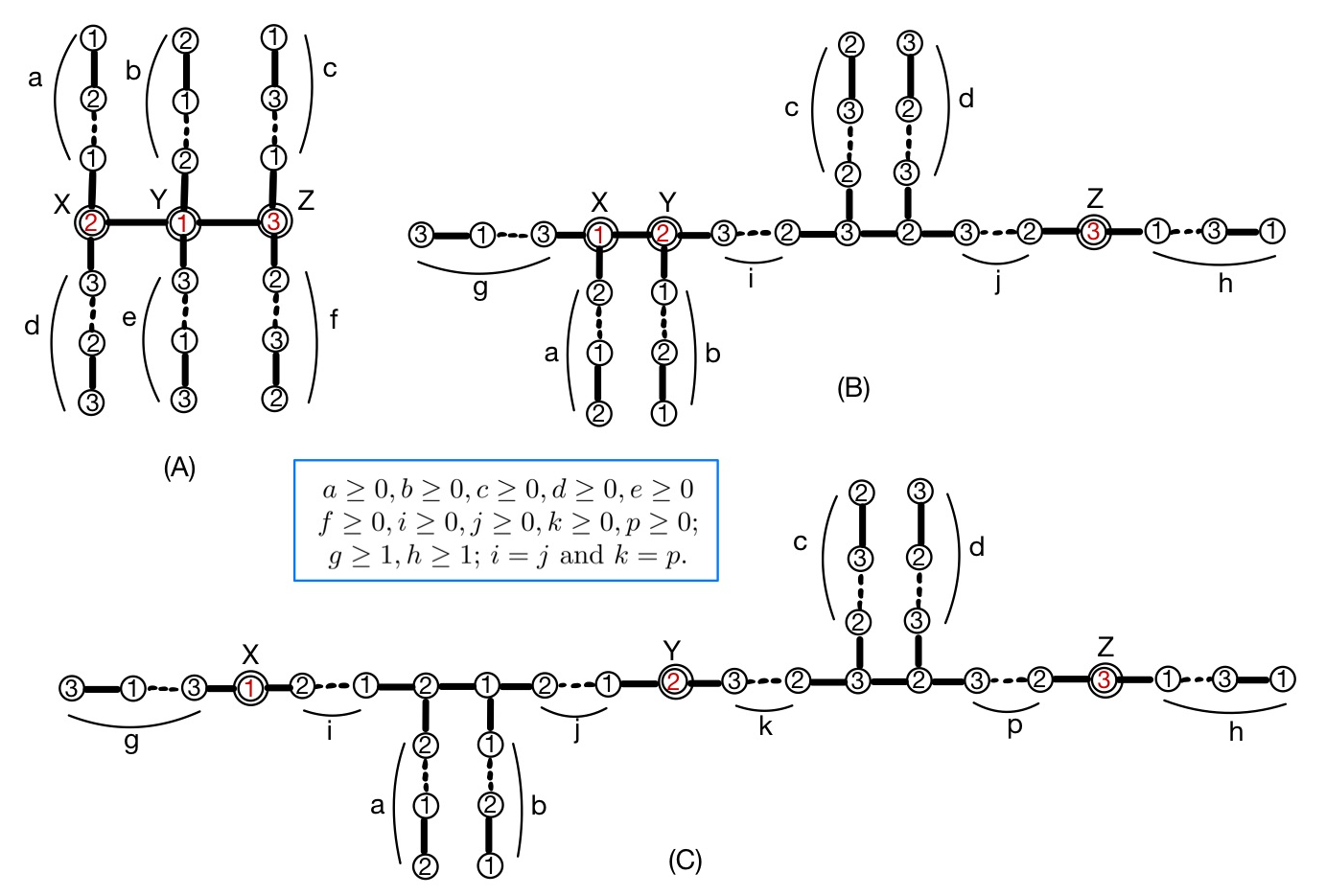}
  \caption{All trees $T$ with $\chi_L(T)=3$ and with a minimum locating coloring.}
\end{center}
\end{figure}

\begin{theorem}
\label{T-tree DL3}
All trees $T \in \mathcal{T}$ with $|Aut(T)| \geq 2$ are the only trees $T$ 
with $\chi_D(T)=\chi_L(T)=3$. 
\end{theorem}

\begin{proof}
Let $T$ be a tree with $\chi_L(T)=3$. Then, $T$ must be isomorphic to one of the trees characterized in Theorem \ref{Th-tree3}. However, not all members $T$ of $\mathcal{T}$ will have $\chi_D(T)=3$. 
If $|Aut(T)| \geq 2$ then any proper 2-coloring of $T$ will be not a distinguishing coloring of $T$. But, by using a locating coloring of such a tree, we have 
$\chi_D(T) = 3$.  If $|Aut(T)|=1$ then any proper 2-coloring of $T$ becomes a distinguishing coloring of $T$ too. So, such a tree $T$ with $|Aut(T)|=1$ will have $\chi_D(T)=2$.
Therefore, we complete the proof.
\end{proof}

Next, we characterize all graphs $G$ other than trees with $\chi_D(G)=\chi_L(G)=3$. Asmiati and Baskoro 
\cite{asm} have characterized all graphs $G$ containing a cycle with $\chi_L(G)=3$. Such graphs are stated in the following theorem.

\begin{theorem}\cite{asm}
\label{Th-graph3}
Let $G$ be a graph with $\chi_L(G)=3$ Then, \\
\begin{enumerate}
\item If $G$ is bipartite then $G$ is isomorphic to any subgraph of the graph in Figure 5 containing at least all blue edges.  
\item If $G$ is not bipartite then $G$ is isomorphic to any subgraph of either the graph (A), (B), (C) or (D) in Figure 6 containing a smallest odd blue cycle $C_m$.\\
\end{enumerate}
\end{theorem}

\vspace{-0.5cm}
\begin{figure}[!h]
\begin{center}
\label{FigE2}
  \includegraphics[width=6cm]{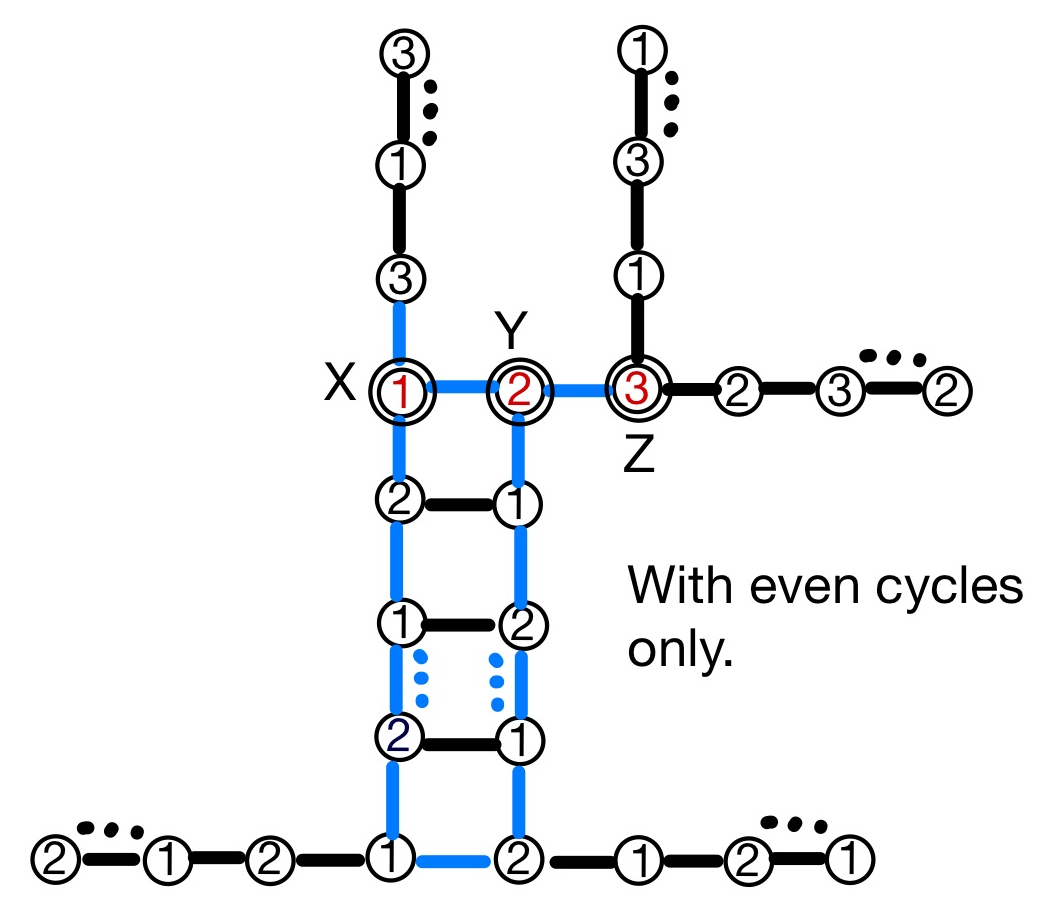}
  \caption{All bipartite graphs $G$ other than trees with $\chi_L(G)=3$ and with a minimum locating coloring.}
\end{center}
\end{figure}

Let $\mathcal{G}$ be the set of all graphs $G$ with $\chi_L(G)=3$, characterized in Theorem \ref{Th-graph3}. Then, we have the following characterization of all graphs $G$ with $\chi_D(G)=\chi_L(G)=3$.

\begin{theorem} Let $G$ be a graph with $\chi_D(G)=\chi_L(G)=3$. Then, $G$ is isomorphic to one of the following:
\begin{enumerate}
\item Any subgraph $H$ of the graph in Figure 4 containing at least all blue edges, with $|Aut(H)| \geq 2$.
\item Any subgraph of either the graph (A), (B), (C) or (D) in Figure 5 
containing a smallest odd blue cycle $C_m$.\\
\end{enumerate}
\end{theorem}

\begin{proof}

Let $G$ be a graph with $\chi_L(G)=3$. Then, $G$ must be isomorphic to one of the graphs characterized in Theorem \ref{Th-graph3}. However, not all members $G$ of $\mathcal{G}$ will have 
$\chi_D(G)=3$. We divide into two cases.\\

\noindent {\bf Case 1}. $G$ is bipartite.\\
If $|Aut(G)| \geq 2$ then any proper 2-coloring of $G$ will be not a distinguishing coloring of $G$. 
But, by using a locating coloring of such a graph, we have $\chi_D(G) = 3$. 
If $|Aut(G)|=1$ then any proper 2-coloring of $G$ becomes a distinguishing coloring of $G$ too. So, such a graph $G$ with $|Aut(G)|=1$ will have $\chi_D(G)=2$. Therefore, all bipartite graphs $G$ in $\mathcal{G}$ with  $|Aut(G)| \geq 2$ will have $\chi_D(G)=3$.\\

\noindent {\bf Case 2}. $G$ is not bipartite.\\
Since $\chi(G) \geq 3$, then by using a locating coloring of $G$ in Figure 6, we obtain $\chi_D(G)=3$. Therefore, all non-bipartite graphs $G$ in $\mathcal{G}$ have $\chi_D(G)=3$. 

Therefore, we complete the proof.
\end{proof}

\begin{figure}[!h]
\begin{center}
\label{FigE3}
  \includegraphics[width=10cm]{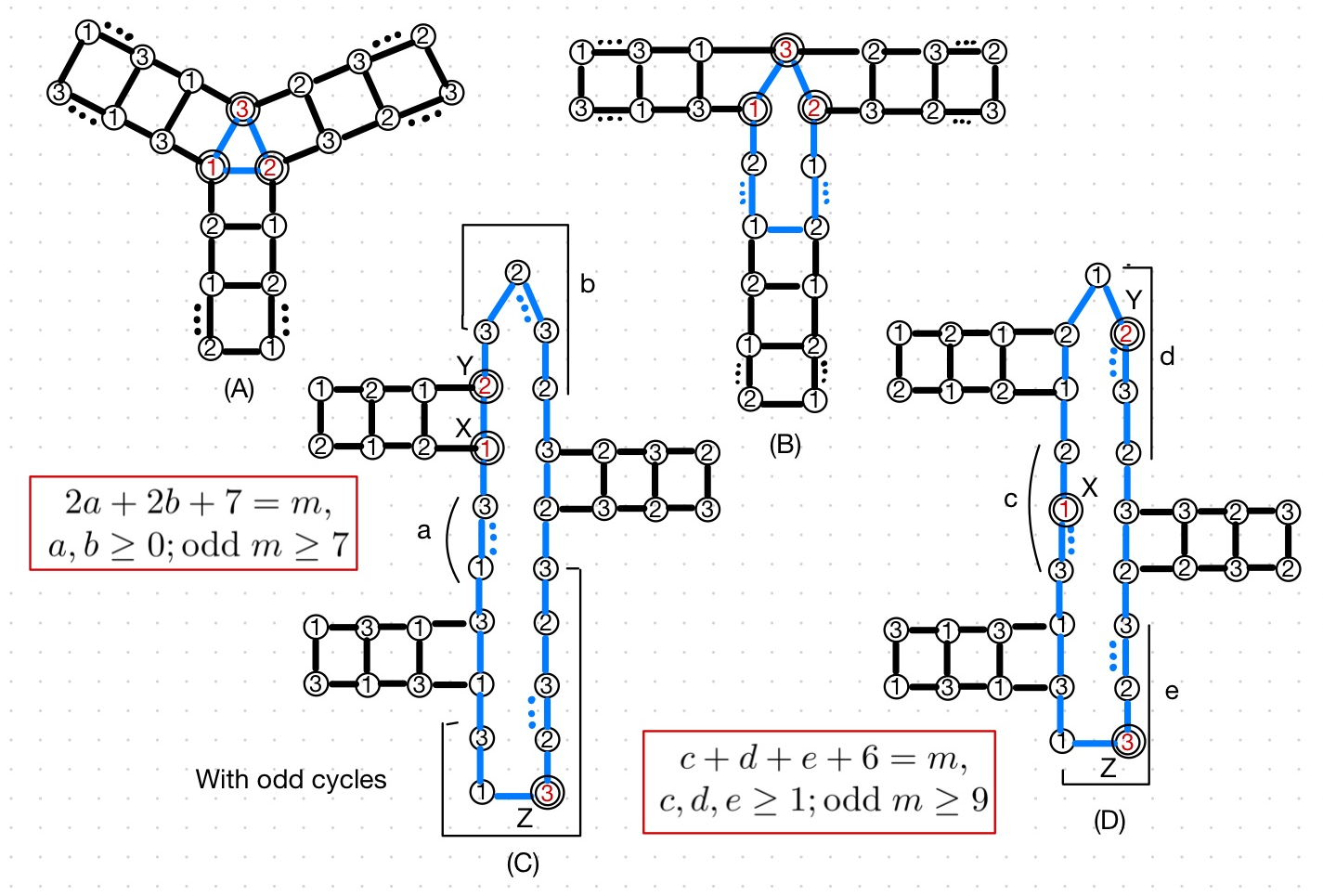}
  \caption{All non-bipartite graphs with $\chi_L(G)=3$ and with a minimum locating coloring.}
\end{center}
\end{figure}

To conclude the paper,  we state some open problem related to all graphs with the same locating and distinguishing chromatic numbers.\\

\noindent {\bf Open Problem.} Characterize a class  {\rsfs A}   of graphs such that 
$ G \in $ {\rsfs A}
if and only if 
$ \chi_{D}(G) = \chi_{L}(G)$. 
\\

\noindent{\bf Acknowledgment.} The authors would like to thank to the anonymous referee for
the valuable comments and suggestions. The third author has been supported by
the World Class Research (WCR) Program, the Indonesian Ministry of Research
and Technology/National Research and Innovation Agency. 

\end{document}